\documentclass[11pt]{article}

\usepackage{amssymb}
\usepackage{amsmath}
\usepackage{amsthm}
\usepackage{graphicx}
\usepackage{siunitx}

\newcommand{\R}{\mathbb R}
\newcommand{\C}{\mathbb C}
\newcommand{\V}{\mathbb V}
\newcommand{\St}{\mathbb S}
\newcommand{\Oh}{{\cal O}}
\newcommand{\st}{{\triangle t}}
\newcommand{\hQ}{{\hat{Q}}}
\newcommand{\hR}{{\hat{R}}}
\renewcommand{\i}{\mathrm i}
\renewcommand{\Im}{\mathrm{Im}}
\renewcommand{\Re}{\mathrm{Re}}
\newcommand{\qedwhite}{\hfill \ensuremath{\Box}}

\newtheorem{remark}{Remark}[section]
\newtheorem{theorem}{Theorem}[section]

\usepackage{xcolor,todonotes}
\author{M.~Schneider, J.~Lang, W.~Hundsdorfer}
\title{Extrapolation-Based Super-Convergent Implicit-Explicit
Peer Methods with A-stable Implicit Part}
\author{
Moritz Schneider\\
{\small \it Technische Universit\"at Darmstadt} \\ {\small \it
Dolivostra{\ss}e 15, 64293 Darmstadt, Germany}\\
{\small moschneider@mathematik.tu-darmstadt.de} \\ \\
Jens Lang\footnote{corresponding author}\\
{\small \it Technische Universit\"at Darmstadt} \\ {\small \it
Dolivostra{\ss}e 15, 64293 Darmstadt, Germany}\\
{\small lang@mathematik.tu-darmstadt.de} \\ \\
Willem Hundsdorfer \\
{\small \it Center for Mathematics and Computer Science} \\
{\small \it P.O.\ Box 94079, 1090 GB Amsterdam, The Netherlands}\\
{\small Willem.Hundsdorfer@cwi.nl}}
\begin{document}
\maketitle

\begin{abstract}
In this paper, we extend the implicit-explicit (IMEX) methods
of Peer type recently developed in [Lang, Hundsdorfer, J. Comp. Phys.,
337:203--215, 2017] to a broader class of two-step methods that allow
the construction of
super-convergent IMEX-Peer methods with A-stable implicit part. IMEX schemes
combine the necessary stability of implicit and low computational costs
of explicit methods to efficiently solve systems of ordinary differential
equations with both stiff and non-stiff parts included in the source term.
To construct super-convergent IMEX-Peer methods with favourable
stability properties, we derive necessary and sufficient conditions on
the coefficient matrices and apply an extrapolation approach based on already computed
stage values. Optimised super-convergent IMEX-Peer methods of
order $s+1$ for $s=2,3,4$ stages are given as result of a search algorithm
carefully designed to balance the size of the stability regions and the extrapolation errors.
Numerical experiments and a comparison to other IMEX-Peer methods
are included.
\end{abstract}


\section{Introduction}
Differential equations of the form $u' = F_0(u) + F_1(u)$, where
$F_0$ is a non-stiff or mildly stiff part and $F_1$ is a stiff
contribution, arise in many initial value problems. Such problems
often result from semi-discretized systems of partial differential
equations with diffusion, advection and reaction terms.
Implicit-explicit (IMEX) methods use this decomposition by treating
only the $F_1$ contribution in an implicit fashion. The advantage
of lower costs for explicit schemes is combined with the favourable
stability properties of implicit schemes to enhance the overall
computational efficiency.

In this paper, we extend the IMEX methods
of Peer type recently developed by Lang and Hundsdorfer
\cite{LangHundsdorfer2017} to a broader class of two-step methods
that include function values from the previous step and thereby allow the
construction of super-convergent IMEX-Peer methods with A-stable implicit part.
Implicit Peer methods have been introduced by Schmitt, Weiner and co-workers
\cite{BeckWeinerPodhaiskySchmitt2012,SchmittWeiner2004,SchmittWeinerErdmann2005}
as a very comprehensive class of general linear methods. Such methods are
described in detail by Butcher \cite{Butcher2006} and Jackiewicz \cite{Jackiewicz2009}. Peer methods are characterized by their special feature that the approximations
in all stages have the same order. They inherit good stability properties and an easy step size change in every time step from one-step methods
without suffering from order reduction for stiff problems. The property that the Peer
stage values have the same order of accuracy can be conveniently exploited to construct
related explicit methods by using extrapolation. The combination of these
implicit and explicit methods leads in a natural way to IMEX methods with the
same order as the original implicit method. This idea was first used by Crouzeix
\cite{Crouzeix1980} with linear multi-step methods of BDF type. Recently,
Cardone, Jackiewicz, Sandu and Zhang \cite{CardoneJackiewiczSanduZhang2014} applied the extrapolation approach to diagonally implicit multistage integration methods and
Lang and Hundsdorfer \cite{LangHundsdorfer2017} to implicit Peer methods constructed
by Beck, Weiner, Podhaisky and Schmitt \cite{BeckWeinerPodhaiskySchmitt2012}.
IMEX-Peer methods are competitive
alternatives to classic IMEX methods for large stiff problems. Higher-order IMEX
Runge-Kutta methods are known to suffer from possible order reduction and serious efficiency loss for stiff problems. A detailed error analysis for increasingly stiff problems has been done by Boscarino \cite{Boscarino2007}, see also the references therein. Moreover, the increasing number of necessary coupling conditions makes their construction difficult.
Kennedy and Carpenter gave general construction principles in the context of additive Runge-Kutta methods \cite{KennedyCarpenter2003}. Boscarino designed an accurate third-order IMEX Runge-Kutta method with a better temporal order of convergence for stiff
problems \cite{Boscarino2009}.

For a method with stage order $q$, it is possible to have convergence
with order equal to $q\!+\!1$. This concept of super-convergence
is discussed for Peer methods applied to non-stiff problems
by Weiner, Schmitt, Podhaisky and Jebens in \cite{WeinerSchmittPodhaiskyJebens2009}
and in the book of Strehmel, Weiner and Podhaisky
\cite[Sect. 5.3]{StrehmelWeinerPodhaisky2012}. It
also holds for stiff problems as shown by Hundsdorfer \cite{Hundsdorfer1994}.
Soleimani and Weiner \cite{SoleimaniWeiner2017a} applied the concept
of super-convergence to an optimally zero-stable subclass of implicit
Peer methods, including also variable time step sizes. More recently,
they have constructed super-convergent IMEX-Peer methods via a partitioning
approach \cite{SoleimaniWeiner2018}. Some of their new methods have
also an A-stable implicit part.
The definition of the order (or order of consistency) of a method commonly
used in Peer literature differs from the more comprehensive definition given
in the book of Hairer, N\o rsett and Wanner
\cite[Sect.\,III.8]{HairerNoersettWanner1993} for general linear methods,
which also covers the super-convergence phenomenon. It uses the concept of
quasi-consistency of Skeel, first introduced in \cite{Skeel1976}. Following this approach,
we will slightly modify the usual criterion for having an extra order
of convergence for Peer methods to construct super-convergent IMEX-Peer
methods based on extrapolation.

The paper is organised as follows. In Section 2, we present the framework to
obtain super-convergent IMEX-Peer methods, combining super-convergent
implicit Peer methods with their corresponding explicit methods derived by
extrapolation. We state necessary and sufficient conditions to ensure the
super-convergence property for this combination. The construction of specific
classes of methods is performed in Section 3. We first summarize all conditions
and then design three new super-convergent IMEX-Peer methods for $s=2,3,4$ with
favourable stability and accuracy properties. Stability regions are given and
compared to those of IMEX-Peer methods from
\cite{LangHundsdorfer2017,SoleimaniKnothWeiner2017}. Numerical results
are presented in Section 4 for a Prothero-Robinson problem, a one-dimensional advection-reaction problem with stiff reactions, and a two-dimensional gravity
wave problem, where selected advection and reaction terms lead to stiffness.

\section{Implicit-Explicit Peer Methods Based on Extrapolation}
\subsection{Super-convergent implicit Peer methods}
We apply the so-called Peer methods introduced by Schmitt, Weiner
and co-workers \cite{SchmittWeiner2004,SchmittWeiner2017,SoleimaniWeiner2017a} to
solve initial value problems in the vector space $\V=\R^m, m\ge 1$,
\begin{equation}
u'(t) = F(u(t)),\quad u(0)=u_0\in\V\,.
\end{equation}
The general form of an $s$-stage implicit Peer method is
\begin{equation}
\label{implpeer1}
w_n = (P\otimes I)w_{n-1} + \st (Q\otimes I)F(w_{n-1}) + \st (R\otimes I)F(w_n)
\end{equation}
with $s\times s$ coefficient matrices $P=(p_{ij})$, $Q=(q_{ij})$,
$R=(r_{ij})$, the $m\times m$ identity matrix $I$, and approximations
\begin{equation}
w_n = [w_{n,1},\ldots,w_{n,s}]^T\in\V^s,\quad w_{n,i}\approx u(t_n+c_i\st)\,,
\end{equation}
where $t_n = n \st$, $n\ge0$, and the nodes $c_i\in\R$ are
such that $c_i\ne c_j$ if $i\ne j$, and $c_s=1$.
Further, $F(w)=[F(w_i)]\in\V^s$ is the application of $F$ to all components
of $w\in\V^s$. The starting
vector $w_0=[w_{0,i}]\in\V^s$ is supposed to be given, or computed
by a Runge-Kutta method, for example.

Peer methods belong to the class of general linear methods. All
approximations have the same order, which gives the name of the methods.
Here, we are interested in A-stable and super-convergent Peer methods with
order of convergence $p\!=\!s+1$, recently constructed by Soleimani and
Weiner in \cite{SoleimaniWeiner2017a}. In the following, for an
$s\times s$ matrix we will use the same symbol for its Kronecker product
with the identity matrix as a mapping from the space $\V^s$ to itself.
Then, (\ref{implpeer1}) simply reads
\begin{equation}
\label{implpeer}
w_n = Pw_{n-1} + \st QF(w_{n-1}) + \st RF(w_n)\,.
\end{equation}
The matrix $R$ is taken to be lower triangular with constant diagonal
$r_{ii}\!=\!\gamma\!>\!0$, $i=1,\ldots,s$, giving singly diagonally
implicit methods. In what follows, we discuss requirements and desirable
properties for the implicit method (\ref{implpeer}).\\

\noindent\textbf{Accuracy.} Let $e=(1,\ldots,1)^T\in\R^s$. We assume
pre-consistency, i.e., $Pe\!=\!e$, which means that for the trivial
equation $u'(t)=0$, we get solutions $w_{n,i}=1$ provided that
$w_{0,j}=1,\;j=1,\ldots,s$. The residual-type local errors result
from inserting exact solution values $w(t_n)=[u(t_n+c_i\st)]\in\V^s$
in the implicit scheme (\ref{implpeer}):
\begin{equation}
r_n = w(t_n) - P w(t_{n-1}) - \st Q w'(t_{n-1}) - \st R w'(t_n)\,.
\end{equation}
Let $c=(c_1,\ldots,c_s)^T$ with
point-wise powers $c^j=(c_1^j,\ldots,c_s^j)^T$. Then Taylor expansion gives
\begin{eqnarray}
w(t_n) &\!=\!& e\otimes u(t_n) + \st c\otimes u'(t_n) +
  \frac{1}{2}\st^2 c^2\otimes u''(t_n) + \ldots\,\\
w(t_{n-1}) &\!=\!& e\otimes u(t_n) + \st (c-e)\otimes u'(t_n) +
  \frac{1}{2}\st^2 (c-e)^2\otimes u''(t_n) + \ldots,
\end{eqnarray}
from which we obtain
\begin{equation}
\label{def:res-r}
r_n = \sum_{j\ge 1} \st^j d_j\otimes u^{(j)}(t_n)
\end{equation}
with
\begin{equation}
\label{def:defect-d}
d_j = \frac{1}{j!}\left( c^j - P(c-e)^j - jQ(c-e)^{j-1} - jRc^{j-1}\right)\,.
\end{equation}
A pre-consistent method is said to have stage order $q$ if
$d_j=0$ for $j=1,2,\ldots,q$. With the Vandermonde matrices
\begin{equation}
\label{eq:V0V1}
V_0 = \big(c_i^{j-1}\big), \qquad
V_1 = \big((c_i-1)^{j-1}\big),\qquad i,j=1,\ldots,s,
\end{equation}
and the diagonal matrices $C=\text{diag}(c_1,c_2,\ldots,c_s)$,
$D=\text{diag}(1,2,\ldots,s)$, the conditions for having stage order
$s$ with the implicit method (\ref{implpeer}) are
\begin{equation}
CV_0 - P(C - I)V_1 - QV_1 D - RV_0 D = 0\,.
\end{equation}
Since $V_1$ and $D$ are regular, we have the relation
\begin{equation}
\label{def:Q}
Q = (CV_0 - P(C - I)V_1 - RV_0 D)(V_1D)^{-1}\,,
\end{equation}
showing that $Q$ is uniquely defined by the choice of $P$, $R$,
and the node vector $c$.\\

\noindent\textbf{Stability.}
Applying the implicit method (\ref{implpeer}) to
Dahlquist's test equation $y'\!=\!\lambda y$ with $\lambda\in\C$,
gives the following recursion for the approximations $w_n$:
\begin{equation}
w_n = (I-zR)^{-1}(P+zQ)w_{n-1}=:M_{im}(z)w_{n-1}
\end{equation}
with $z\!:=\!\lambda\st$. Hence, $w_n=M_{im}(z)^nw_{0}$. The
matrix $P\!=\!M_{im}(0)$ should be power bounded to have stability
for the trivial equation $u'(t)\!=\!0$. This requirement of
zero-stability is enforced by Schmitt, Weiner et al. by taking
$P$ such that one eigenvalue equals $1$ (due to pre-consistency)
and the others are $0$. Such methods are called optimally zero-stable.
We will also look at methods that are A-stable, i.e., the spectral radius of $M_{im}(z)$ satisfies $\rho(M_{im}(z))\le 1$ for all $z\in\C$ with $\Re (z)\le 0$. Since $M_{im}(\infty)\!=\!R^{-1}Q$ with $Q\neq 0$, A-stability does not imply
L-stability. To guarantee good damping properties for very stiff problems, we
will aim at having a small spectral radius of $R^{-1}Q$.\\

\noindent\textbf{Super-convergence.}
We are interested in using the degrees of freedom provided by
the free parameters in $P$, $R$, and $c$ to have convergence
with order $p\!=\!s+1$
without raising the stage order further. This is discussed under
the heading super-convergence in the book of Strehmel, Weiner and
Podhaisky \cite[Sect.\,5.3]{StrehmelWeinerPodhaisky2012}
for non-stiff problems. It is related to the definition of order of
consistency for general linear methods as given in
\cite[Sect.\,III.8]{HairerNoersettWanner1993}.
Similar results for stiff systems can be found in \cite{Hundsdorfer1994}.
According to the last paper, we will slightly modify the usual criterion for
having an extra order of convergence for Peer methods to later construct
super-convergent IMEX-Peer methods based on extrapolation.

Let $\varepsilon_n\!=\!w(t_n)-w_n$ be the global error. Under the
standard stability assumption, where products of the transfer matrices
are bounded in norm by a fixed constant $K$ (see, e.g., Theorem 2 in \cite{SoleimaniWeiner2017a}), we get the estimate
$\|\varepsilon_n\|\le K(\|\varepsilon_0\|+\|r_1\|+\ldots+\|r_n\|)$.
Together with stage order $s$, this gives the standard convergence result
\begin{equation}
\label{implpeer-globerr}
\|\varepsilon_n\| \le K\|\varepsilon_0\|+\st^s t_nK\|d_{s+1}\|_\infty
\max_{0\le t\le t_n}\|u^{(s+1)}(t)\| + \Oh(\st^{s+1})\,.
\end{equation}
Then we have the following
\begin{theorem}
\label{Th:super-implpeer}
Assume the implicit Peer method (\ref{implpeer}) has stage order $s$ and estimate (\ref{implpeer-globerr}) holds true for the global error.
Then the method is convergent of order $p\!=\!s$. Furthermore,
if $d_{s+1}\in\text{range}\,(I-P)$ and the initial values are of order
$s+1$, then the order of convergence is $p\!=\!s+1$.
\end{theorem}
\noindent {\it Proof}: The first statement follows directly from
(\ref{implpeer-globerr}). Suppose that $d_{s+1}=(I-P)v$ with $v\in\R^s$, and let
\begin{equation}
\bar{w}(t_n) := w(t_n) - \st^{s+1}v\otimes u^{(s+1)}(t_n)\,.
\end{equation}
Insertion of these modified solution values in the scheme (\ref{implpeer}) will give modified local errors
\begin{equation}
\begin{array}{rll}
\bar{r}_n &=& \bar{w}(t_n)-P\bar{w}(t_{n-1})-\st QF(\bar{w}(t_{n-1}))
-\st R F(\bar{w}(t_n))\\[2mm]
&=& r_n - \st^{s+1}d_{s+1}\otimes u^{(s+1)}(t_n)+\Oh(\st^{s+2})\,,
\end{array}
\end{equation}
which, due to (\ref{def:res-r}), reveals $\bar{r}_n\!=\!\Oh(\st^{s+2})$. For $\bar{\varepsilon}_n=\bar{w}(t_n)-w_n$ this yields, in the same way as above,
$\|\bar{\varepsilon}_n\|\le K\|\varepsilon_0\|+\Oh(\st^{s+1})$.
Since $\|\bar{\varepsilon}_n-\varepsilon_n\|\le \st^{s+1}\|v\|_\infty\|u^{(s+1)}(t_n)\|$ and $\|\varepsilon_0\|=\Oh(\st^{s+1})$, this shows convergence of order $s\!+\!1$ for the global errors $\varepsilon_n$. \qedwhite\\

Recall that the range of $I-P$ consists of the vectors that are orthogonal to the null space of $I-P^T$. If the method is zero-stable, then this null space has dimension one. So up to a constant there is a unique vector $v\in\R^s$ such that $(I-P^T)v=0$.
Then we have
\begin{equation}
\label{super-cond}
d_{s+1}\in\text{range}\,(I-P)\quad \text{iff}\quad v^Td_{s+1}\!=\!0\,.
\end{equation}
We fix $v$ by $v^Te\!=\!1$ and set $P\!=\!ev^T$ to ensure
pre-consistency ($Pe=e$), optimal zero-stability and
$(I-P^T)v=0$. In this way, $P$ is
determined by the vector $v$, which has to satisfy the conditions
\begin{equation}
v^Te=1\quad\text{and}\quad v^Td_{s+1}=0
\end{equation}
to achieve super-convergence of order $s\!+\!1$.

A closer inspection of the global error $\varepsilon_n$ reveals
that the condition $P^jd_{s+1}=0$ for all $j\ge s-1$ is also
an appropriate way to construct super-convergent Peer methods.
This approach has been used by Schmitt, Weiner et al.
\cite{SchmittWeiner2017,WeinerSchmittPodhaiskyJebens2009} for
explicit and implicit schemes. There is a strong relation to
(\ref{super-cond}), which can be stated in the following
\begin{theorem}
\label{Th:ppower0}
(1) Let $P^jd_{s+1}=0$ for all $j\ge s-1$ and
$eig(P)=\{1,\lambda_2,\ldots,\lambda_s\}$
with $|\lambda_i|<1, i=2,\ldots,s$, i.e., the Peer method is
zero-stable. Then, $d_{s+1}$ lies in the range of $I-P$.
(2) Suppose $d_{s+1}$ lies in the range of $I-P$ and
$eig(P)=\{1,0,\ldots,0\}$, i.e., the Peer method is
optimally zero-stable. Then $P^jd_{s+1}=0$ for all $j\ge s-1$.
\end{theorem}
\noindent {\it Proof}: Due to the spectra of $P$ in both cases, there is a
regular matrix $S$ such that
\begin{equation}
\label{p_jordan}
P = S^{-1}
\begin{pmatrix}
1 & 0\\ 0 & J
\end{pmatrix}
S =: S^{-1}\hat{P}S\,,
\end{equation}
where $J\in\R^{(s-1)\times (s-1)}$ consists of certain Jordan blocks.
(1) Using this form, the assumption on powers of $P$ reads
$S^{-1}\hat{P}^jSd_{s+1}=0$
for $j\ge s-1$, from which follows that the first component of
the vector $Sd_{s+1}$ vanishes. Then, the solvability of the
equation $(I-P)x=S^{-1}(I-\hat{P})Sx=d_{s+1}$ reduces to the
question whether $I_{s-1}-J$ with $I_{s-1}$ being the unit matrix in
$\R^{(s-1)\times (s-1)}$ is invertible. This is indeed the case,
since it is an upper triangular matrix with values $1-\lambda_i\neq 0$,
$i=2,\ldots,s$ on the diagonal. (2) In this case, the assumption gives
the existence of an $x\in\R^s$ such that $(I-P)x=d_{s+1}$. Using
(\ref{p_jordan}), we have $(I-\hat{P})Sx=Sd_{s+1}$, which shows that
the first component of the vector $Sd_{s+1}$ vanishes. Since $I-\hat{P}$
is a strictly upper diagonal matrix, we deduce
$P^jd_{s+1}=S^{-1}\text{diag}(1,0,\ldots,0) Sd_{s+1}=0$ for
$j\ge s-1$.\qedwhite

\begin{remark}
The second statement in Theorem \ref{Th:ppower0} does not hold, if the
Peer method is zero-stable, but not optimally zero-stable. In this case,
there exist super-convergent Peer methods which do not satisfy the
conditions $P^jd_{s+1}=0$ for all $j\ge s-1$. An example is the new
IMEX-Peer4s method constructed in Section \ref{Sec:Construction}.
\end{remark}

\subsection{Extrapolation}
In \cite{LangHundsdorfer2017}, Lang and Hundsdorfer have applied extrapolation to an
implicit method of the form (\ref{implpeer}) with $Q\!=\!0$ and order $s$
to find a related explicit method and eventually construct
IMEX-Peer methods of order $s$ with good stability
properties. This procedure is well-known from linear multistep
methods, see for instance Crouzeix \cite{Crouzeix1980} or the review
in the book of Hundsdorfer and Verwer
\cite[Sect.\,IV.4.2]{HundsdorferVerwer2003}. It was also used
by Cardone, Jackiewicz, Sandu, and Zhang \cite{CardoneJackiewiczSanduZhang2014}
to construct implicit-explicit diagonally implicit multistage integration methods.
Here, we will extend this extrapolation idea to obtain super-convergent
IMEX-Peer methods of order $s\!+\!1$, where the implicit method is
A-stable and the stability region of the overall method is optimised.

Starting with an implicit method (\ref{implpeer}), where all approximations
$w_{n,j}$ have order $s$, we can obtain a corresponding explicit method by
extrapolation of $F(w_n)$ with order $s$. Using $w_{n-1}$ and most recent values $w_{n,j},j=1,\ldots,i-1$, available in the $i$th-stage with $1<i<s$,
gives
\begin{equation}
\label{expol_f}
F(w_n) = S_1F(w_{n-1})+S_2F(w_n)+\Oh(\st^{s})\,,
\end{equation}
where the extrapolation coefficients are collected in the
$s\!\times\!s$ matrices $S_1$ and $S_2$. Note that $S_2$ is strictly
lower triangular. Applied in (\ref{implpeer}), this leads to an explicit
method of the form
\begin{equation}
\label{explpeer}
w_n = Pw_{n-1} + \st\hQ F(w_{n-1}) + \st\hR F(w_n)\,,
\end{equation}
with $\hQ=Q+RS_1$ and the strictly lower triangular matrix $\hR=RS_2$,
since $R$ is lower triangular. In what follows, we discuss properties and
the issue of super-convergence for this explicit method.\\

\noindent\textbf{Accuracy.} With exact values $F(w(t_n))\in\V^s$, the
residual-type error vector for the extrapolation can be expanded by
Taylor series at $t_n$:
\begin{equation}
\begin{array}{rll}
\delta_n &=& F(w(t_n))-S_1F(w(t_{n-1}))-S_2F(w(t_n)) \\[2mm]
 &=& \displaystyle\sum_{j\ge 0}\frac{1}{j!}\left( (I-S_2)c^j-S_1(c-e)^j
\right)\otimes\frac{d^j}{dt^j}F(u(t_n))\,\st^j\,.
\end{array}
\end{equation}
The conditions for order $s$ read
\begin{equation}
\label{eq:stage_order_s}
(I-S_2) c^j - S_1 (c-e)^j = 0,\quad 0\le j\le s-1\,,
\end{equation}
which can be reformulated to $S_1=(I-S_2)V_0V_1^{-1}$. Thus, the choice of a strictly
lower triangular $S_2$ determines $S_1$.

Using the expression for $\delta_n$, the conditions for stage order $s$ of the implicit method, and (\ref{eq:stage_order_s}), we derive for the residual-type local error of the explicit method (\ref{explpeer}) the form
\begin{equation}
r_n = \st^{s+1}\left( d_{s+1}+R \,l_{s}\right)\otimes u^{(s+1)}(t_n) +\Oh(\st^{s+2})
\end{equation}
with
\begin{equation}
\label{def:expol-ls}
l_{s} = \frac{1}{s!}\left( (I-S_2)c^{s}-S_1(c-e)^s\right)\,.
\end{equation}
Thus, by construction, all the stages have again order $s$, at least,
so (\ref{explpeer}) is an explicit Peer method.\\

\noindent\textbf{Super-convergence.} First, we proceed as for the implicit
method under the standard stability assumption. With stage order $s$ of
the implicit method (\ref{implpeer}) and order $s$ of the extrapolation
in (\ref{expol_f}), we analogously get the convergence result for the global
error of the explicit method defined in (\ref{explpeer}),
\begin{equation}
\label{explpeer-globerr}
\|\varepsilon_n\| \le K\|\varepsilon_0\|+\st^s t_nK\|d_{s+1}+R\,l_s\|_\infty
\max_{0\le t\le t_n}\|u^{(s+1)}(t)\| + \Oh(\st^{s+1})\,.
\end{equation}
Then we have the following
\begin{theorem}
Assume the implicit Peer method (\ref{implpeer}) has stage order $s$,
conditions (\ref{eq:stage_order_s}) for the extrapolation are satisfied,
and estimate (\ref{explpeer-globerr}) holds true for the global error.
Then the explicit Peer method (\ref{explpeer}) is convergent of order $p\!=\!s$. Furthermore, if $(d_{s+1}+R\,l_s)\in\text{range}\,(I-P)$ and the initial values are of order $s+1$, then the order of convergence is $p\!=\!s+1$.
\end{theorem}
\noindent {\it Proof}: The first statement follows directly from
(\ref{explpeer-globerr}). Replacing $d_{s+1}$ by $d_{s+1}+R\,l_s$
in the proof of Theorem \ref{Th:super-implpeer} gives the desired
result. \qedwhite\\

With this result, we can conclude, in the same way as above for the implicit
method, that super-convergence is achieved if
\begin{equation}
v^T (d_{s+1} + R\,l_s)=0 \quad\text{with }v\in\R^s\text{ such that } (I-P^T)v=0\,.
\end{equation}
If the implicit Peer method is already super-convergent, this simplifies to
$v^TR\,l_s\!=\!0$.

\subsection{Super-convergent IMEX-Peer methods}
Combining the related implicit and explicit methods (\ref{implpeer}) and
(\ref{explpeer}) yields an IMEX method for systems of the form
\begin{equation}
\label{ode}
u'(t) = F_0(u(t)) + F_1(u(t))\,,
\end{equation}
where $F_0$ will represent the non-stiff or mildly stiff part, and $F_1$ gives
the stiff part of the equation. The resulting IMEX scheme is
\begin{equation}
\label{imex-peer}
w_n = P w_{n-1} + \st\hQ F_0(w_{n-1}) + \st\hR F_0(w_{n})
+ \st QF_1(w_{n-1}) + \st R F_1(w_{n}),
\end{equation}
where the extrapolation is used only on $F_0$.
Here, $\hQ=Q+RS_1$ and $\hR=RS_2$. For non-stiff problems,
this IMEX method will have order $s$ for any decomposition $F=F_0+F_1$. However,
for stiff problems it should be required that the derivatives of
$\varphi_k(t)=F_k(u(t)), \;k=0,1,$ are bounded by a moderate constant which is
not affected by the stiffness parameters, such as the spatial mesh width $h$ for
semi-discrete systems obtained from PDEs.

With exact solution values $u(t_{n,i})$, standard local consistency analysis
for the IMEX-Peer method (\ref{imex-peer}) gives for the residual-type
local errors
\begin{equation}
\label{imexpeer-res}
r_n = E_{im} + \st R\,E_{ex} + \Oh\left(\st^{s+2}\right)\,,
\end{equation}
with
\begin{equation}
E_{im} = \st^{s+1} d_{s+1}\otimes u^{(s+1)}(t_n)\quad\text{and}
\quad E_{ex} = \st^s l_s\otimes\frac{d^s}{dt^s}F_0(u(t_n))\,.
\end{equation}
Therefore, by standard convergence arguments, we have the following
\begin{theorem}
\label{Th:imex-conv}
Let the $s$-stage implicit Peer method (\ref{implpeer}) with
coefficients $(c,P,Q,R)$, $Q$ defined in (\ref{def:Q}), be zero-stable
and suppose its stage order is equal to $s$.
Let the initial values satisfy $w_{0,i}-u(t_0+c_i\st)=\Oh(\st^s),\,
i=1,\ldots,s$.
Then the IMEX-Peer method (\ref{imex-peer}) with $\hR=RS_2$ and
$\hQ=Q+R(I-S_2)V_0V_1^{-1}$ is convergent of order $s$ for constant step size
and arbitrary strictly lower triangular matrix $S_2$.
\end{theorem}
Combining the requirements for super-convergence of order $s+1$ stated
above for the explicit and implicit Peer methods we have
\begin{theorem}
\label{Th:imex-superconv}
Let the assumptions of Theorem \ref{Th:imex-conv} be fulfilled and
the IMEX-Peer method (\ref{imex-peer}) be convergent of order $s$.
If the initial values are of order $s+1$, $d_{s+1}\in\text{range}\,(I-P)$
and $R\,l_{s}\in\text{range}\,(I-P)$, then the order of convergence
is $s+1$.
\end{theorem}
\noindent {\it Proof}: Suppose $d_{s+1}=(I-P)v_d$ and $R\,l_s=(I-P)v_l$
with $v_d,v_l\in\R^s$, and let
\begin{equation}
\bar{w}(t_n) = w(t_n) - \st^{s+1}v_d\otimes u^{(s+1)}(t_n) -
\st^{s+1}v_l\otimes \frac{d^s}{dt^s}F_0(u(t_n))\,.
\end{equation}
Inserting these values in (\ref{imex-peer}) gives the modified
residual-type local errors
\begin{equation}
\begin{array}{rll}
\bar{r}_n &=& \bar{w}(t_n) - P \bar{w}(t_{n-1}) - \st\hQ F_0(\bar{w}(t_{n-1}))
- \st\hR F_0(\bar{w}(t_n))\\[2mm]
 && - \st QF_1(\bar{w}(t_{n-1})) - \st R F_1(\bar{w}(t_n))\,,
\end{array}
\end{equation}
which can be rearranged to
\begin{equation}
\begin{array}{rll}
\bar{r}_n &=& \bar{w}(t_n) - P \bar{w}(t_{n-1}) - \st Q F(\bar{w}(t_{n-1}))
- \st R F(\bar{w}(t_n))\\[2mm]
&& + \st R \left( F_0(\bar{w}(t_n)) - S_1F_0(\bar{w}(t_{n-1}))
- S_2F_0(\bar{w}(t_n)) \right)\,.
\end{array}
\end{equation}
Taylor expansions gives
\begin{equation}
\bar{r}_n =
\displaystyle r_n - \st^{s+1}d_{s+1}\otimes u^{(s+1)}(t_n)
- \st^{s+1}R\,l_s\otimes \frac{d^s}{dt^s}F_0(u(t_n)) + \Oh(\st^{s+2})
\end{equation}
with $r_n$ as defined in (\ref{imexpeer-res}). This shows
$\bar{r}_n=\Oh(\st^{s+2})$. Then, the same arguments as in the proof
of Theorem \ref{Th:super-implpeer} give convergence of order $s+1$ for
the global errors $\varepsilon_n=w(t_n)-w_n$.\qedwhite\\

Let $d_{s+1}\in\text{range}\,(I-P)$ and $R\,l_{s}\in\text{range}\,(I-P)$.
Then, with the unique vector $v\in\R^s$ such that $(I-P^T)v=0$ and $v^Te=1$
it holds
\begin{equation}
y\in\text{range}\,(I-P)\quad\text{iff}\quad v^Ty=0
\end{equation}
with $y=d_{s+1},\,R\,l_{s}$. Setting $P=ev^T$, we enforce
pre-consistency, optimal zero-stability and $(I-P^T)v=0$.
Furthermore, $P$ is fully determined by
the vector $v$. We will use this simplifying construction to find
suitable super-convergent IMEX-Peer methods for $s=2,3$. To
enrich the space of suitable matrices $P$ for $s=4$,
we only request zero-stability and follow the
approach discussed in Theorem {\ref{Th:imex-superconv}
with $d_{s+1},R\,l_s\in\text{range}(I-P)$ to achieve
super-convergence.

\subsection{Stability of IMEX-Peer methods}
In order to discuss stability properties of the IMEX-Peer method
(\ref{imex-peer}), we consider the split scalar test equation
\begin{equation}
\label{imex-test-eq}
y'(t) = \lambda_0 y(t) + \lambda_1y(t),\quad t\ge 0,
\end{equation}
with complex parameters $\lambda_0$ and $\lambda_1$.
Applying an IMEX-Peer method to (\ref{imex-test-eq}) gives the recursion
\begin{equation}
\label{imex-peer-lin}
w_{n+1} = (I-z_0RS_2-z_1R)^{-1}(P+z_0Q+z_0RS_1+z_1Q)w_n
=:M(z_0,z_1)w_n
\end{equation}
with $z_i=h\lambda_i, i=0,1$. Therefore, stability is ensured if
\begin{equation}
\label{imex-peer-stabmat}
\rho (M(z_0,z_1)) < 1.
\end{equation}
The stability regions of the IMEX-Peer method for
$\alpha \in [0^\circ,90^\circ]$ are defined by the sets
\begin{equation}
\St_\alpha = \{z_0\in\C: (\ref{imex-peer-stabmat})\text{ holds for any }
z_1 \in\C \text{ with }|\Im (z_1)|\le -\tan(\alpha)\cdot\Re(z_1) \}
\end{equation}
in the left-half complex plane. Further, we define the stability region
of the corresponding explicit method as
\begin{equation}
\St_E = \{z_0\in\C: \rho (M(z_0,0)) < 1\}
\end{equation}
with the stability matrix $M(z_0,0)=(I-z_0RS_2)^{-1}(P+z_0Q+z_0RS_1)$.
Efficient numerical algorithms to compute $\St_\alpha$ and $\St_E$ are
extensively described
in \cite{CardoneJackiewiczSanduZhang2014,LangHundsdorfer2017}.

Since $\St_\alpha\subset \St_E$, the goal is to construct IMEX-Peer methods
for which $\St_E$ is large and $\St_E\backslash \St_\alpha$ is as small as possible for angles $\alpha$ that are close to $90^\circ$. In what follows, we will use
the additional degrees of freedom provided by introducing the terms
$\st QF(w_{n-1})$ in the implicit method to construct super-convergent
IMEX-Peer methods with a non-empty stability region $\St_{90^\circ}$, i.e., the
underlying implicit Peer method is A-stable.

\section{Construction of Super-Convergent IMEX-Peer Methods Based on Extrapolation}
\label{Sec:Construction}
We will first summarize the conditions on the coefficients of the methods derived in
the previous sections and give a procedure for the construction of
super-convergent IMEX-Peer methods with favourable stability properties
for $s=2,3,4$.\\

\begin{table}[h!]
\centering
\begin{tabular}{|l r c|c l r|}
\hline
\multicolumn{6}{|l|}{IMEX-Peer2s, $s=2$, optimally zero-stable}\\
\hline
$c_1$	    & $0.591977499693304$ 	&&& $p_{11}$ & 	$-1.082167419515352$ \\
$c_2$		& $1.000000000000000$ 	&&&	$p_{12}$ & 	 $2.082167419515352$ \\	
$\gamma$	& $0.969486340522434$ 	&&& $p_{21}$ & 	$-1.082167419515352$ \\
$r_{21}$	& $-1.007885680522306$	&&& $p_{22}$ & 	 $2.082167419515352$ \\
$s_{21}$    & $0.819167640511257$   &&& &		\\
\hline\hline							
\multicolumn{6}{|l|}{IMEX-Peer3s, $s=3$, optimally zero-stable}\\
\hline
$c_1$		&   $0.173922498101250$ 	&&&		$p_{11}$	&$-0.516269158723393$ \\
$c_2$		&   $0.584759944717930$ 	&&&		$p_{12}$	& $2.301256858880021$  \\
$c_3$		&   $1.000000000000000$ 	&&&		$p_{13}$	&$-0.784987700156628$ \\
$\gamma$	&   $0.456150901216430$     &&&		$p_{21}$	&$-0.516269158723393$ \\
$r_{21}$	&   $0.271188675194957$     &&&		$p_{22}$	& $2.301256858880021$  \\
$r_{31}$	&   $0.099808771568803$     &&&		$p_{23}$	&$-0.784987700156628$ \\
$r_{32}$	&   $0.395734854902157$     &&&		$p_{31}$	&$-0.516269158723393$  \\
$s_{21}$    &   $1.500000000000000$     &&&		$p_{32}$	& $2.301256858880021$  \\
$s_{31}$    &   $0.204731875658678$     &&&		$p_{33}$ 	&$-0.784987700156628$ \\		$s_{32}$    &   $1.320000000000000$     &&& &   \\
\hline\hline
\multicolumn{6}{|l|}{IMEX-Peer4s, $s=4$, zero-stable}\\
\hline
$c_1$		&  $-0.926697334544583$ 	&&&    	$p_{11}$ 	&  $0.164346920652337$ \\
$c_2$		&   $0.180751924024702$ 	&&&		$p_{12}$	&  $1.941408294648193$ \\
$c_3$		&   $0.850343633101352$     &&&		$p_{13}$	& $-2.764059964877189$ \\
$c_4$		&   $1.000000000000000$     &&&		$p_{14}$	&  $1.658304749576660$ \\
$\gamma$	&   $0.413154106969917$     &&&     $p_{21}$	&  $0.424734281438207$ \\
$r_{21}$	&   $1.186201415903827$     &&&     $p_{22}$	&  $1.133423589655944$ \\
$r_{31}$	&   $1.327861645060559$     &&&     $p_{23}$	& $-0.792340606563880$ \\
$r_{32}$	&   $0.525143168803633$     &&&     $p_{24}$	&  $0.234182735469729$ \\
$r_{41}$	&   $1.324984727912657$     &&&     $p_{31}$	&  $0.562642125818718$ \\
$r_{42}$	&   $0.576558985833141$     &&&     $p_{32}$	&  $0.131525283967289$ \\
$r_{43}$	&   $0.071014878172581$     &&&     $p_{33}$	&  $2.162128869126546$ \\
$s_{21}$	&   $3.884803988586850$     &&&     $p_{34}$	& $-1.856296278912553$ \\
$s_{31}$	&  $-3.053336552626494$     &&&     $p_{41}$	&  $0.589388877693458$ \\
$s_{32}$	&   $2.821635541838257$     &&&     $p_{42}$	& $-0.169092459871472$ \\
$s_{41}$	&  $-3.555025951383727$     &&&     $p_{43}$	&  $3.071031564759426$ \\
$s_{42}$	&   $2.895140468767150$     &&&     $p_{44}$	& $-2.491327982581412$\\
$s_{43}$	&   $0.162040780709875$     &&&  & \\
\hline
\end{tabular}
\parbox{13cm}{
\caption{Coefficients of the super-convergent $s$-stage IMEX-Peer methods IMEX-Peer2s,
IMEX-Peer3s, and IMEX-Peer4s for $s=2,3,4$ with $S_2=(s_{ij})$.}
\label{tab:coeffs}
}
\end{table}

\noindent\textbf{Implicit method.} An
$s$-stage implicit Peer method is determined by the coefficient matrices
$P,Q,R\in\R^{s\times s}$, and the node vector $c\in\R^s$. We look for
singly diagonally implicit methods, i.e., $R$ is taken to be lower triangular
with a constant $\gamma>0$ on the diagonal. Stage order $s$ is imposed through,
see (\ref{def:Q}),
\begin{equation}
Q = (CV_0 - P(C - I)V_1 - RV_0 D)(V_1D)^{-1}\,.
\end{equation}
The matrix $P$ is chosen such that the method is pre-consistent, (optimally)
zero-stable, and super-convergent. More precisely, we set
\begin{itemize}
\item[(1)] $P=ev^T$ with $v\in\R^s$ and $v^Te=1$ for $s=2,3$,
\item[(2)] $P$ such that the method is zero-stable for $s=4$.
\end{itemize}
In both cases, super-convergence is obtained by satisfying $v^Td_{s+1}=0$
with $d_{s+1}$ defined in (\ref{def:defect-d}) and $v$ such that
$(I-P^T)v\!=\!0$. Note
that the special choice of $P$ in (1) yields an optimally zero-stable method,
whereas for case (2) the property of zero-stability has to be incorporated in the
search algorithm. We were not able to find suitable methods within the setting
(1) for $s=4$. Further, for the nodes $c_i\in\R$, we use $c_i\in (0,1]$ for
$s=2,3$ and $c_i\in(-1,1]$ for $s=4$.\\

\noindent\textbf{Explicit method.} The extrapolation is uniquely defined by
the choice of the matrices $S_1$ and $S_2$ in (\ref{expol_f}). Order $s$ is
guaranteed by setting, see (\ref{eq:stage_order_s}),
\begin{equation}
S_1 = (I-S_2)V_0V_1^{-1}\,.
\end{equation}
So the entries of $S_2\in\R^{s\times s}$ are free parameters. Super-convergence
is obtained by satisfying $v^TR\,l_s=0$ with
$l_s$ defined in (\ref{def:expol-ls}) and $v$ such that $(I-P^T)v\!=\!0$.

\begin{table}[ht]
\centering
\begin{tabular}{|l|crrc|rr|cc|}
\hline\rule{0mm}{5mm}\hspace{-0.1cm}
IMEX- & $\alpha$ & $|\St_\alpha|$ & $x_{max}$ & $\rho(R^{-1}Q)$ & $|\St_E|$ & $y_{max}$ & $c_{im}$ & $c_{ex}$\\[1mm]
\hline\rule{0mm}{5mm}\hspace{-0.1cm}
Peer2s & $90.0^\circ$ & $2.15$ & $-1.41$ & $1.28\,10^{-1}$ & $4.47$ & $1.21$ & $2.37\,10^{-1}$ & $3.23\,10^{-1}$\\[2mm]
Peer3s & $90.0^\circ$ & $2.67$ & $-1.58$ & $5.52\,10^{-1}$ & $6.11$ & $1.69$ & $1.24\,10^{-1}$ & $1.68\,10^{-1}$\\[2mm]
Peer4s & $90.0^\circ$ & $1.07$ & $-1.45$ & $5.42\,10^{-1}$ & $4.39$ & $1.00$ & $6.42\,10^{-2}$ & $1.17\,10^{-1}$\\[1mm]
\hline\rule{0mm}{5mm}\hspace{-0.1cm}
Peer2  & $90.0^\circ$ &  $7.44$ & $-4.86$ & $0.00$ &  $8.53$ & $0.40$ & $7.05\,10^{-2}$ & $2.78\,10^{-1}$\\[2mm]
Peer3  & $86.1^\circ$ &  $8.28$ & $-3.07$ & $0.00$ & $10.68$ & $1.78$ & $8.20\,10^{-3}$ & $3.58\,10^{-2}$\\[2mm]
Peer4  & $83.2^\circ$ &  $4.64$ & $-3.57$ & $0.00$ & $9.36$ & $1.90$ & $3.43\,10^{-4}$ & $4.27\,10^{-3}$\\[1mm]
\hline\rule{0mm}{5mm}\hspace{-0.1cm}
Peer3a  & $90.0^\circ$ &  $2.87$ & $-1.69$ & $1.60\,10^{-3}$ & $5.19$ & $1.75$ & $1.46\,10^{-1}$ & $1.90\,10^{-1}$\\[2mm]
Peer4a  & $90.0^\circ$ &  $2.65$ & $-1.73$ & $1.24\,10^{-1}$ & $3.53$ & $1.15$ & $8.97\,10^{-2}$ & $1.41\,10^{-1}$\\[1mm]
\hline
\end{tabular}\\
\parbox{13cm}{
\caption{\small
Size of stability regions $\St_\alpha$ and $\St_E$, $x_{max}$ at the negative
real axis, $y_{max}$ at the positive imaginary axis, spectral radius of $R^{-1}Q$,
and error constants $c_{im}=\|d_{s+1}\|$ and $c_{ex}=\|R\,l_s\|$ for IMEX-Peer methods with $s=2,3,4$, including those from
\cite{LangHundsdorfer2017,SoleimaniKnothWeiner2017}.}
\label{tab-stab-imex}
}
\end{table}

\noindent\textbf{Construction principles.} We start with the search for a
super-convergent implicit Peer method along the following design criteria:\\[2mm]
\phantom{a}\hspace{1cm} A-stability, $\rho(R^{-1}Q)$ is close to zero,
$\|P\|_F,\|Q\|_F,\|R\|_F,\|d_{s+1}\|$ are small.\\[2mm]
This is done using the {\sc Matlab}-routine {\it fminsearch}, where
we include the desired properties in the objective function and use random
start values for the remaining degrees of freedom. Different combinations
of weights in the objective function are employed to select promising
candidates which are then used for a subsequent extrapolation. The latter
aims at finding parameters in $S_2$ such that we have the following
properties:
\\[2mm]
\phantom{a}\hspace{1cm} large stability regions $\St_E$ and $\St_{90^\circ}$,
$\|S_1\|_F,\|S_2\|_F,\|R\,l_s\|$ are small.\\[2mm]
Again {\it fminsearch} is used with different combinations of weights in
the objective function. Following this approach, we have found new
super-convergent IMEX-Peer methods for $s=2,3,4$. The coefficients of
the methods for
$c,P,R$ and $S_2=(s_{ij})$ are given in Table~\ref{tab:coeffs}.

The resulting values for the stability regions $\St_\alpha$ and $\St_E$ as well
as other constants are collected in Table~\ref{tab-stab-imex}. For
comparison, we also show the values for the IMEX-Peer methods tested
in \cite{LangHundsdorfer2017} (IMEX-Peer2, IMEX-Peer3, IMEX-Peer4) and
in \cite{SoleimaniKnothWeiner2017} (IMEX-Peer3a, IMEX-Peer4a), where
IMEX-Peer-3a is nearly super-convergent in the implicit part. It can be
observed that super-convergence comes with (i) smaller stability regions
and (ii) significantly larger error constants for higher order. However,
the convergence is one order higher and A-stability of the
implicit method pays off for problems with eigenvalues on the imaginary
axis as can be seen from our numerical experiments.
More details on the stability regions are shown in
Figure~\ref{fig-stabreg-all-shape}.

\begin{figure}[t]
\setlength{\unitlength}{1cm}
\centering
\includegraphics[width=0.32\textwidth]{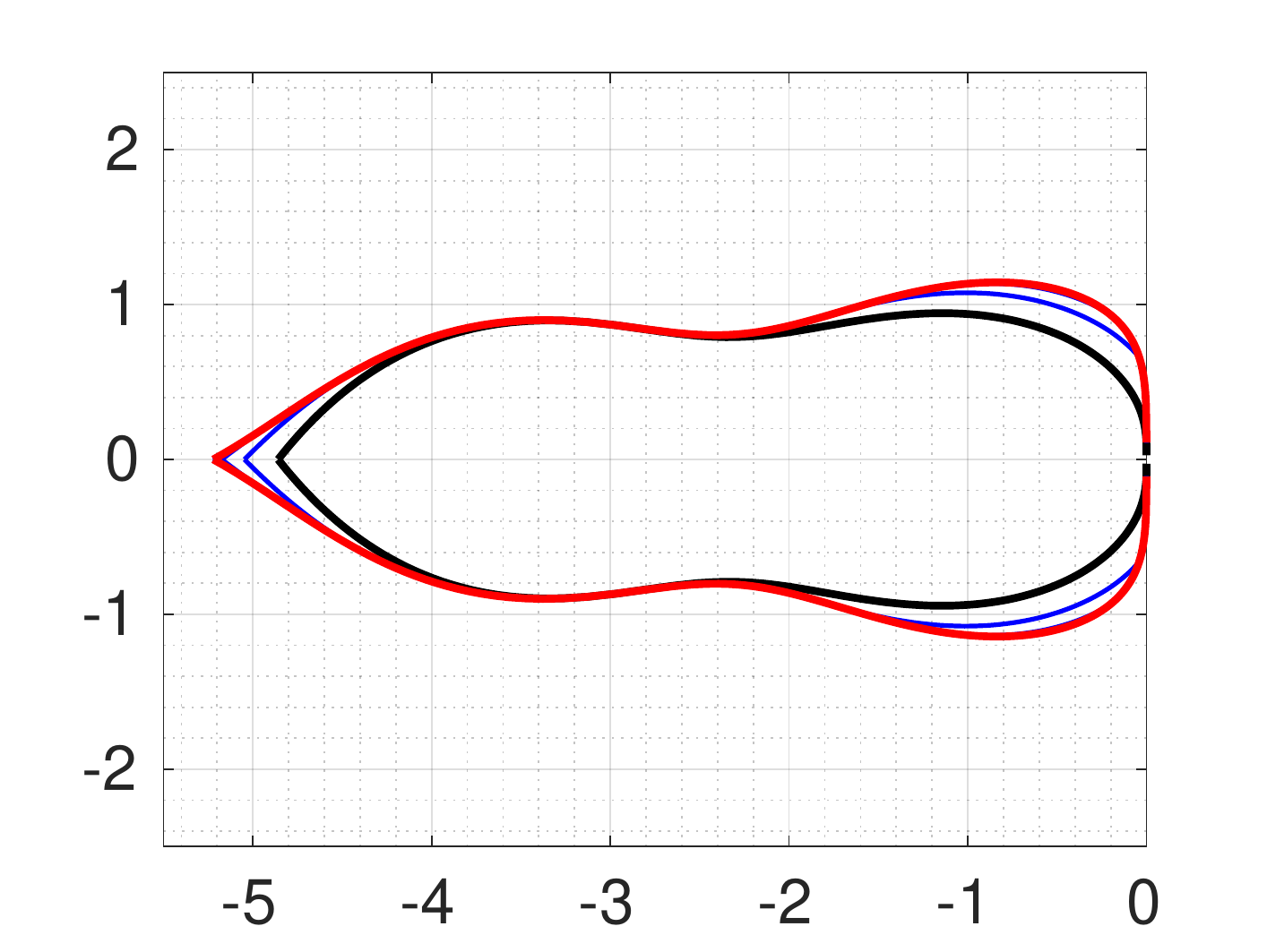}
\put(-3.8,2.8){\scriptsize IMEX-PEER2}
\includegraphics[width=0.32\textwidth]{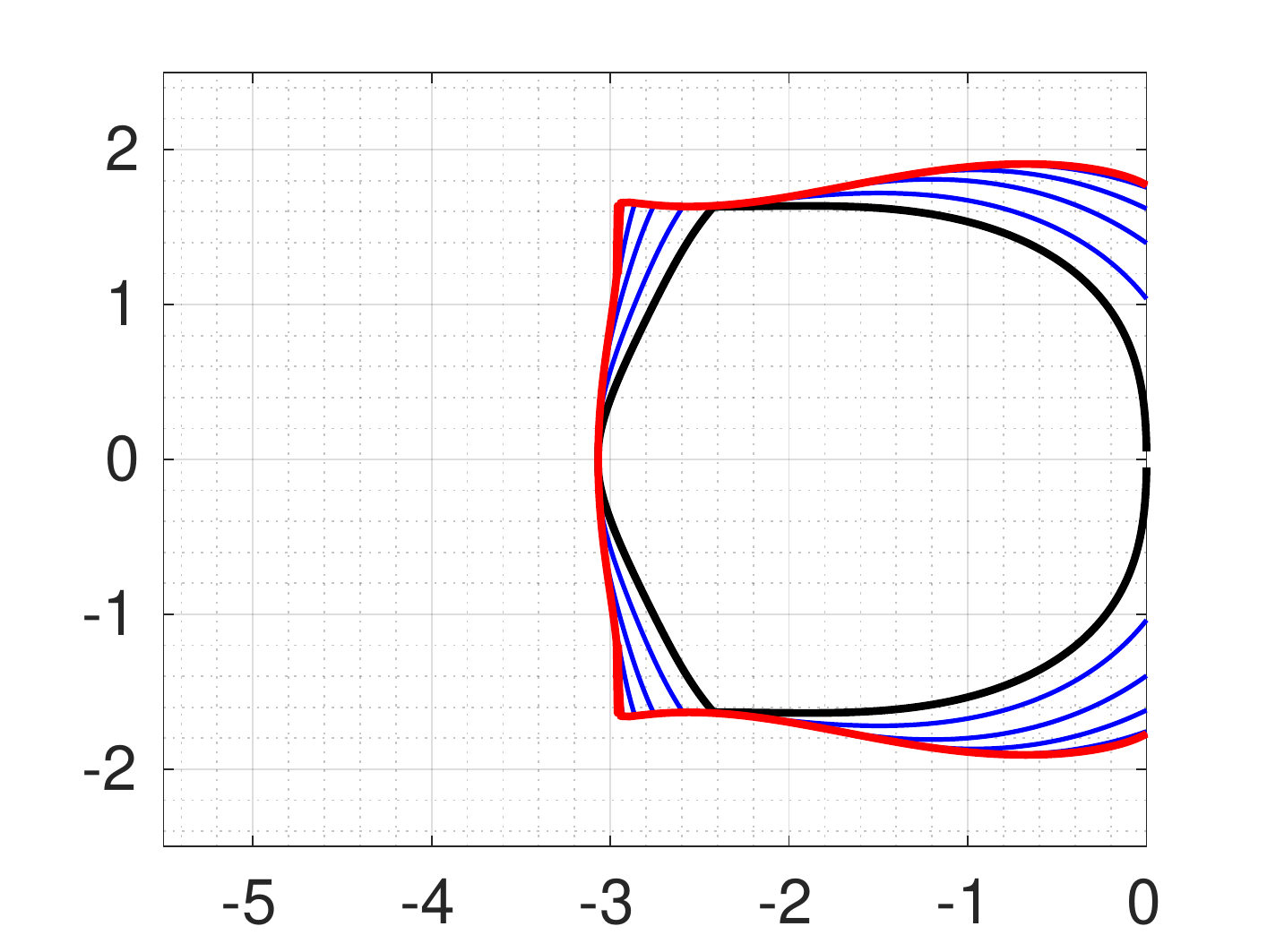}
\put(-3.8,2.8){\scriptsize IMEX-PEER3}
\includegraphics[width=0.32\textwidth]{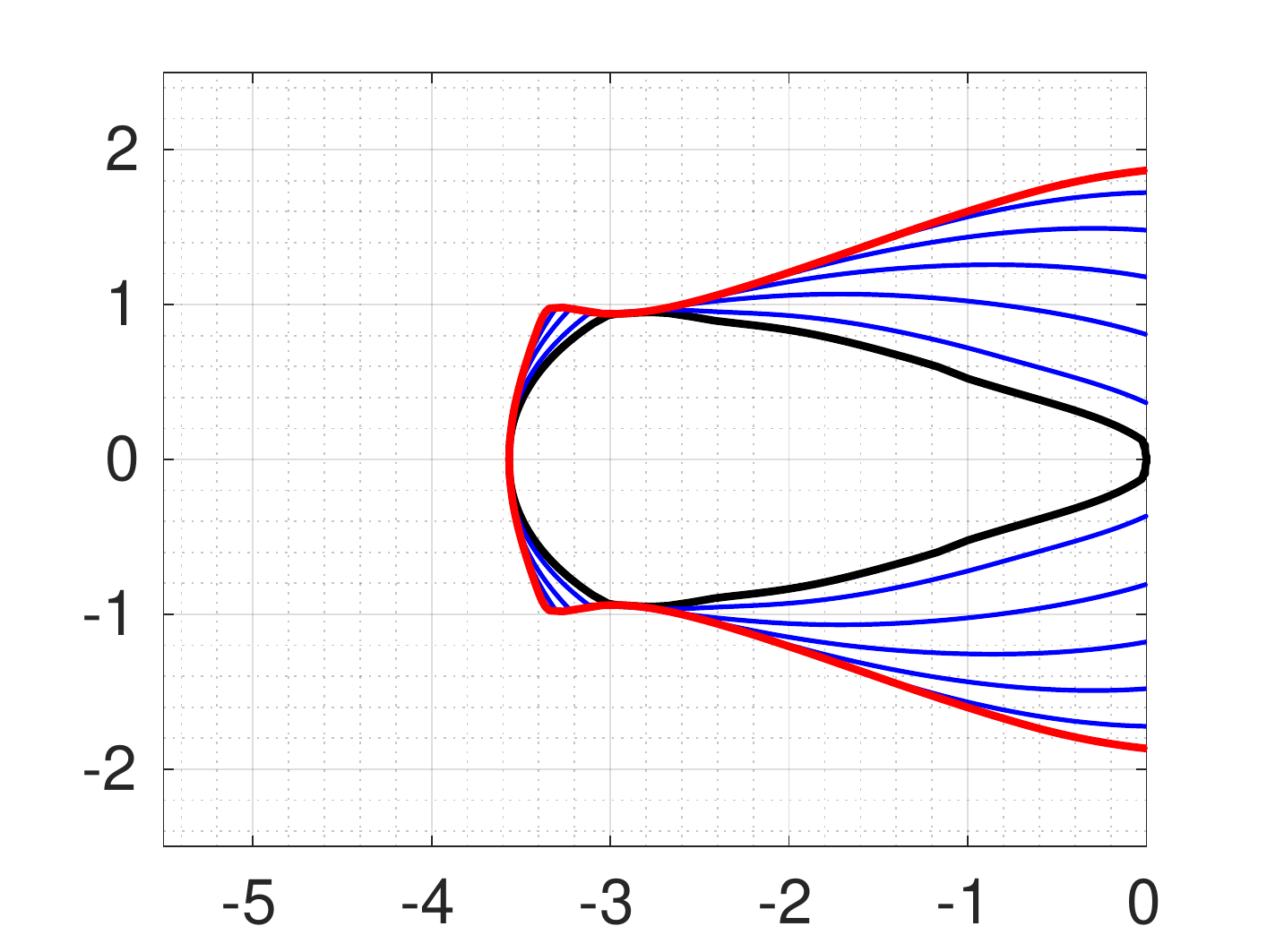}
\put(-3.8,2.8){\scriptsize IMEX-PEER4}
\\
\includegraphics[width=0.32\textwidth]{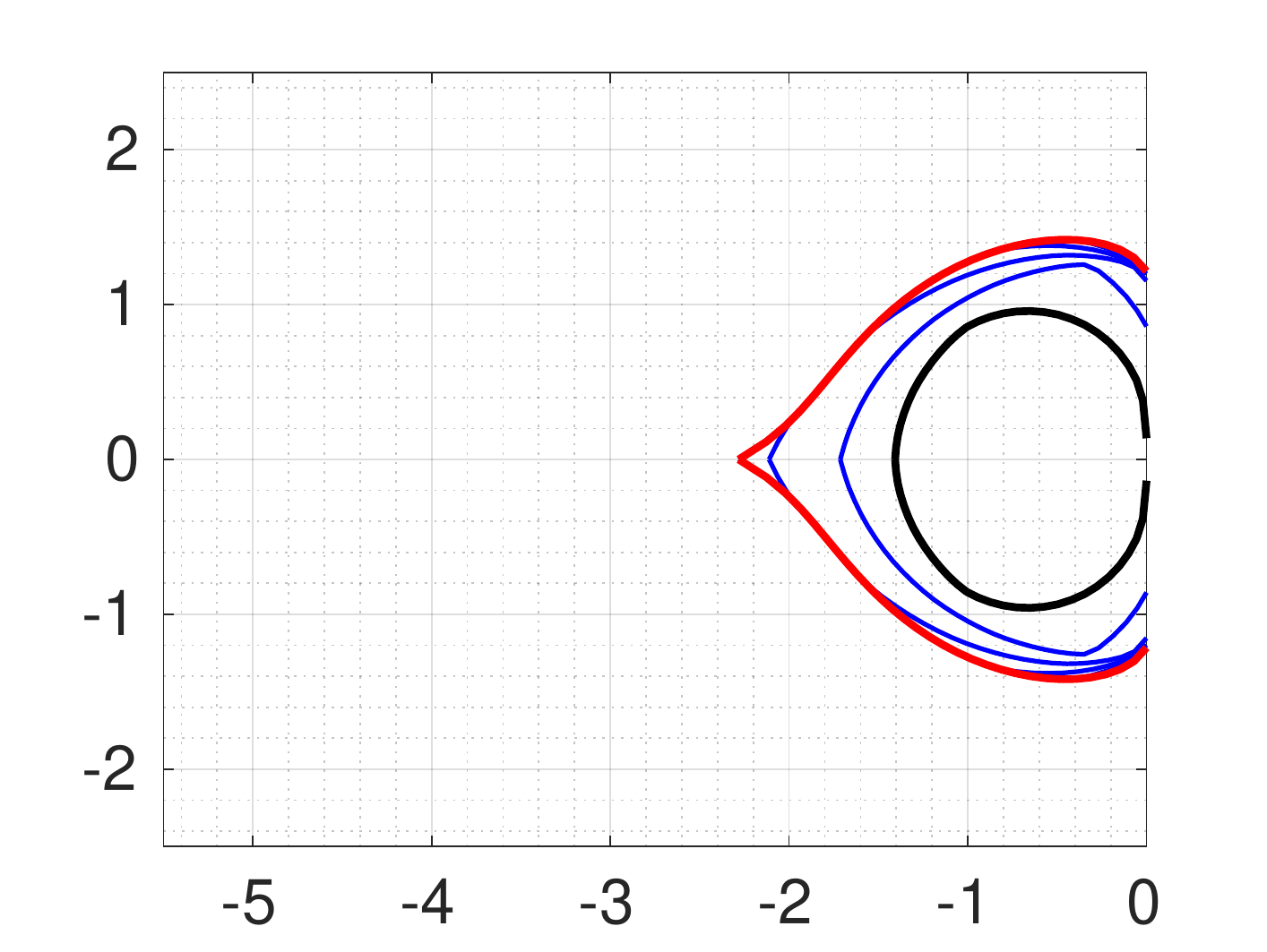}
\put(-3.8,2.8){\scriptsize IMEX-PEER2s}
\includegraphics[width=0.32\textwidth]{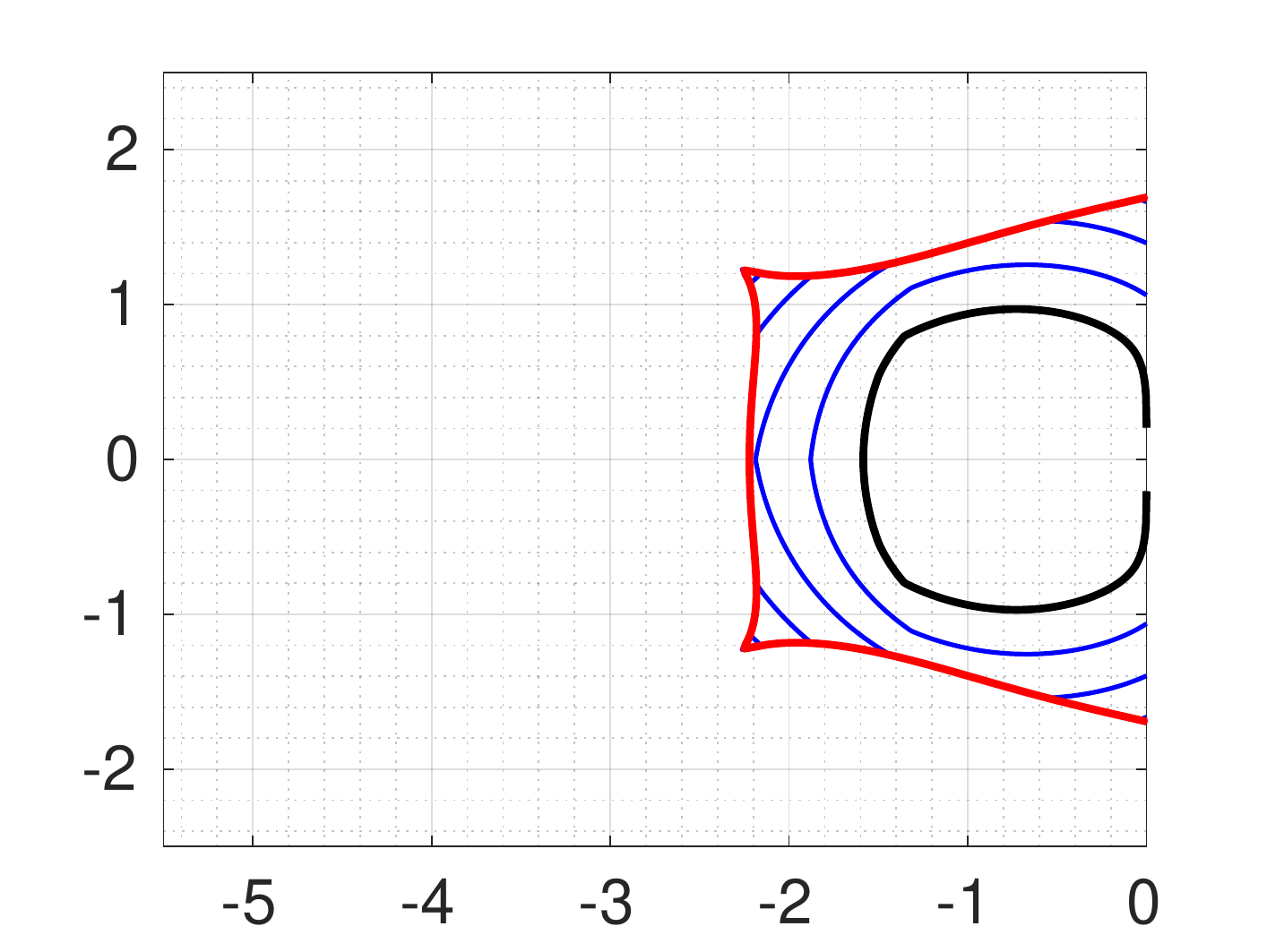}
\put(-3.8,2.8){\scriptsize IMEX-PEER3s}
\includegraphics[width=0.32\textwidth]{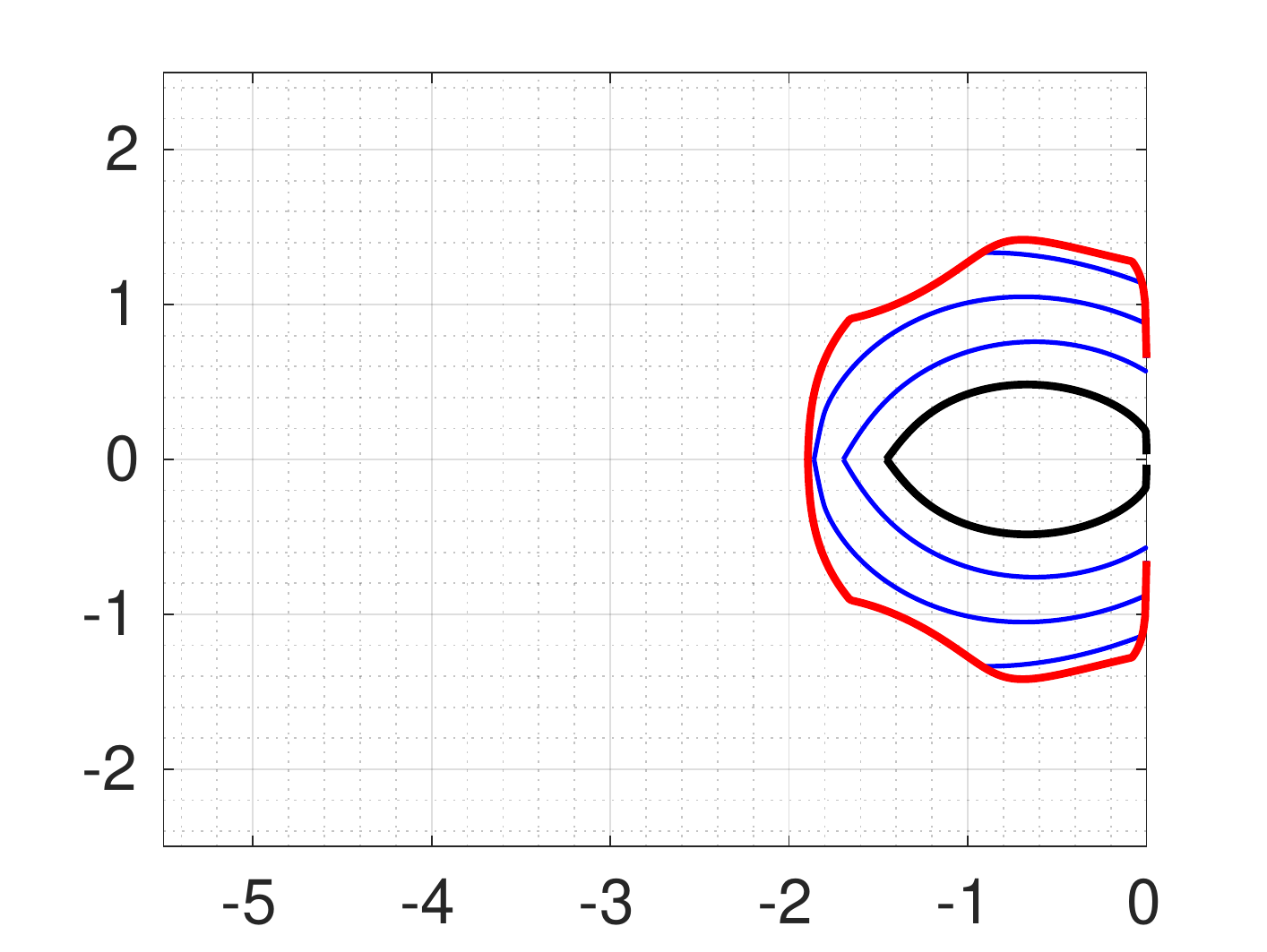}
\put(-3.8,2.8){\scriptsize IMEX-PEER4s}
\\
\hspace{4.3cm}
\includegraphics[width=0.32\textwidth]{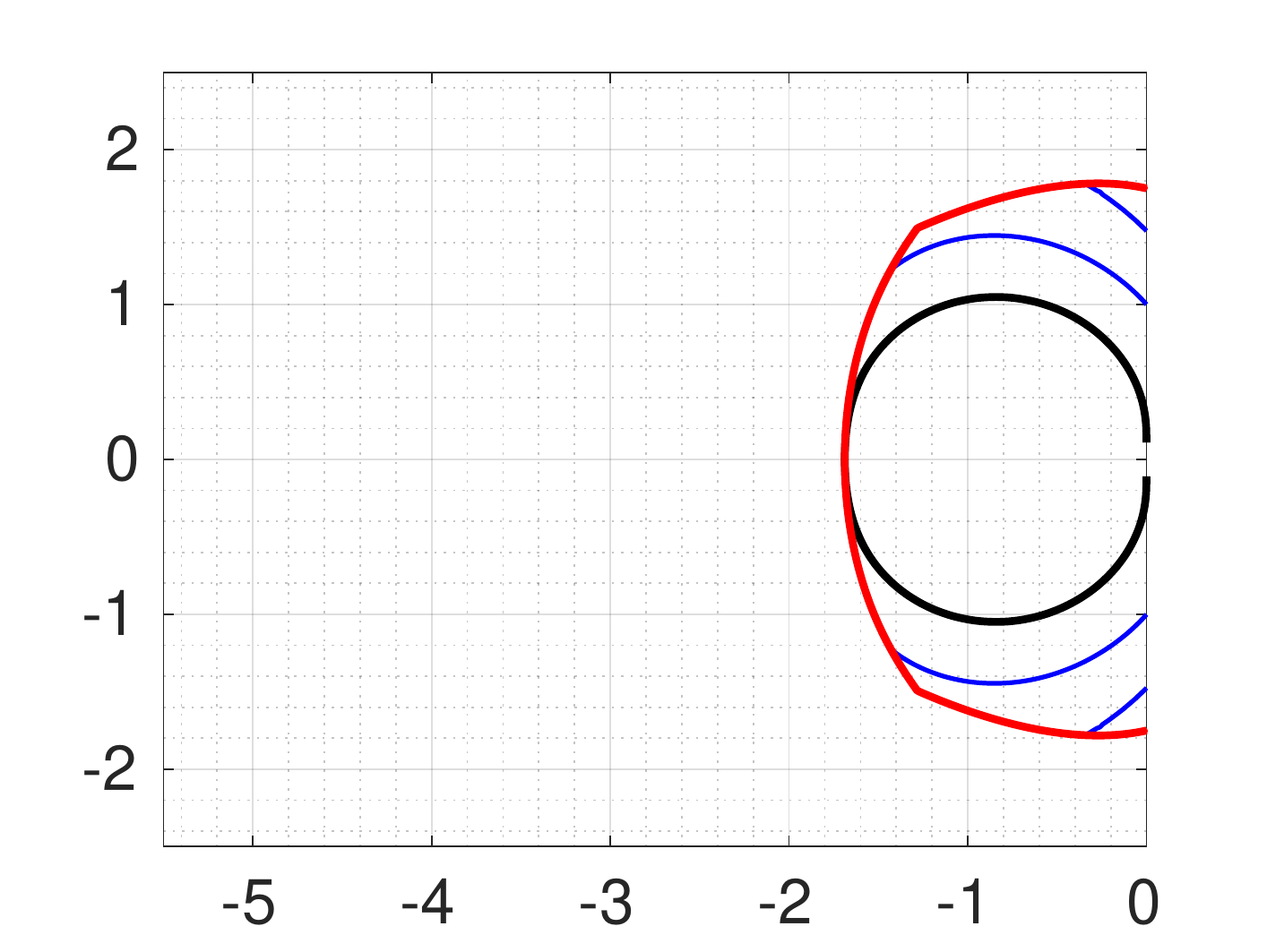}
\put(-3.8,2.8){\scriptsize IMEX-PEER3a}
\includegraphics[width=0.32\textwidth]{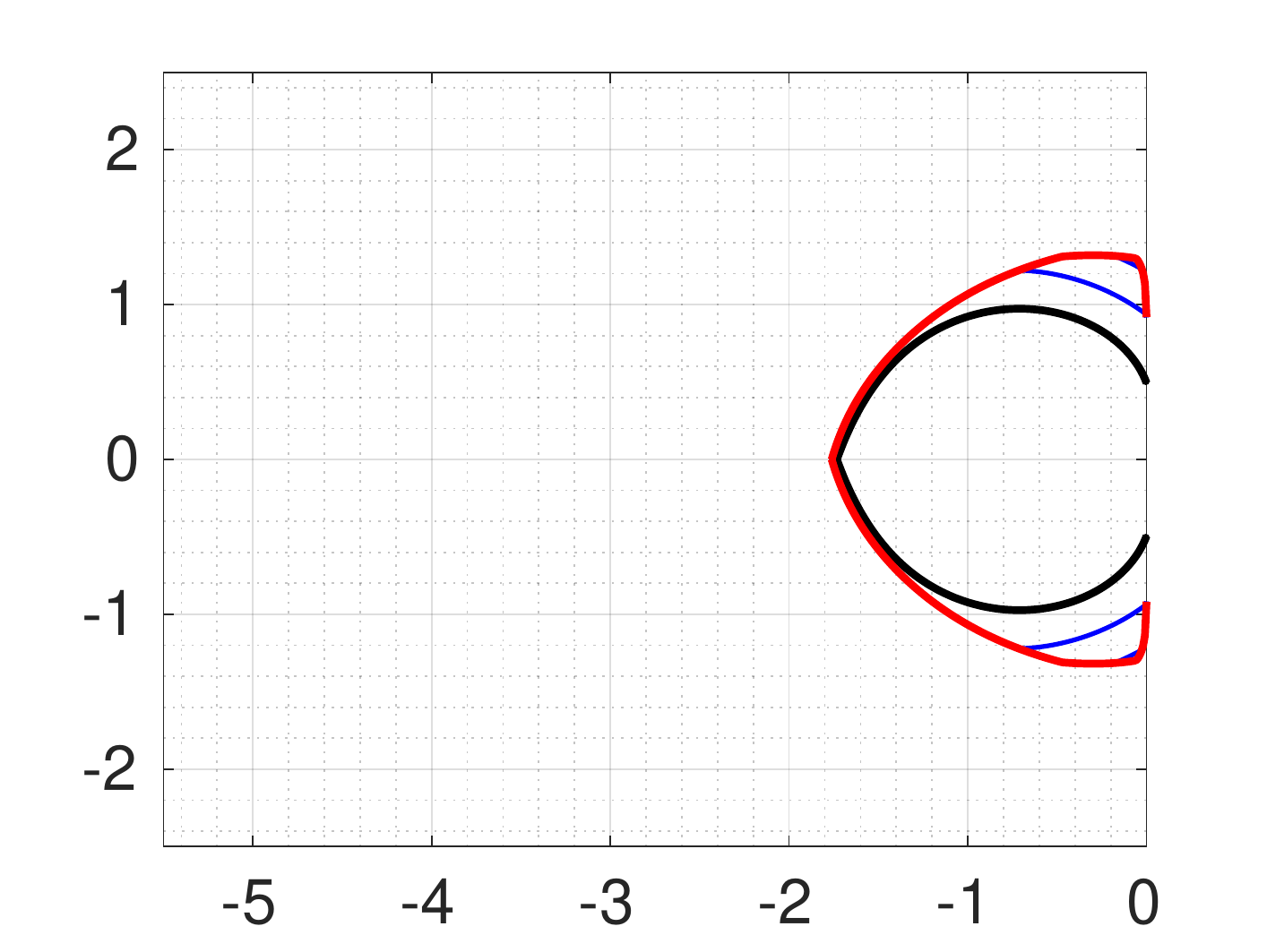}
\put(-3.8,2.8){\scriptsize IMEX-PEER4a}
\\
\parbox{13cm}{
\caption{Stability regions
$\St_{\alpha}$ (black line), $\St_{\beta}$ for
$\beta=75^\circ,60^\circ,45^\circ,30^\circ,15^\circ$ (blue lines),
and $\St_0$ (red line) for IMEX-Peer methods with $s=2,3,4$ (left to right).
We have $\alpha=90^\circ$, i.e., A-stability of the underlying implicit method,
for all methods except IMEX-Peer3 ($\alpha=86.1^\circ$)
and IMEX-Peer4 ($\alpha=83.3^\circ$).}
\label{fig-stabreg-all-shape}
}
\end{figure}

\section{Numerical Experiments}

\subsection{Prothero-Robinson Problem}
In order to study the rate of convergence under stiffness, we consider the
Prothero-Robinson type equation used in \cite{SoleimaniKnothWeiner2017},
\begin{align}
y' & =
\begin{pmatrix}
0 \\ y_1 + y_2 - \sin(t) \end{pmatrix} + \begin{pmatrix}
-10^6 (y_1 - \cos(t)) + 10^3 (y_2 -\sin(t)) - \sin(t) \\ 0
\end{pmatrix}\,,
\end{align}
where $t\in [0,5].$ The first term is treated explicitly and the second
implicitly. Initial values are taken from the analytic solution
$y(t)\!=\!(\cos(t),\sin(t))^T$. The error of the approximate solution
$Y$ is calculated at the final time $T\!=\!5$ in the scaled
maximum norm, i.e., $err=\max_{i=1,2}|Y_i-y_i|/(1+|y_i|)$. The
values for $\st=5/(100+60i),i=0,\ldots,8$ are shown in Figure~\ref{fig:ProRob}.

All methods show their theoretical order of convergence quite nicely.
The smaller error constants of the IMEX-Peer(s) methods
for $s=3,4$ compared to the super-convergent IMEX-Peer methods with
the same order are also visible.

\begin{figure}
\centering
\includegraphics[width=0.8\linewidth]{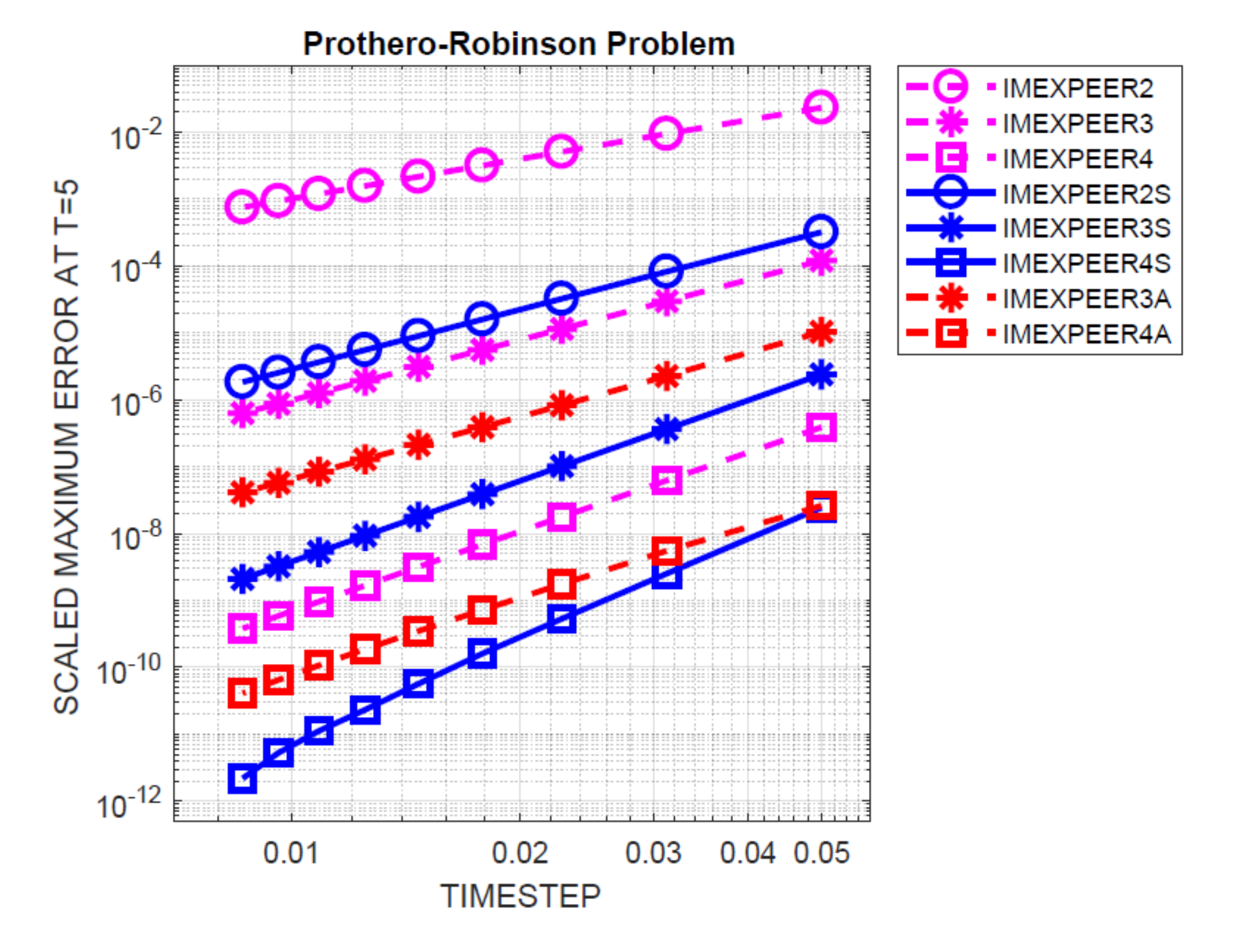}
\parbox{13cm}{
\caption{
Prothero-Robinson Problem: Scaled maximum errors at $T=5$
vs. time step sizes. Comparison of IMEX-Peer methods. The
convergence orders derived from a least squares fit are:
$1.95, 2.94, 2.99, 3.14, 4.00, 3.94, 3.68, 5.21$ (top to bottom).}
\label{fig:ProRob}
}
\end{figure}

\subsection{Linear Advection-Reaction Problem}
A second PDE problem for an accuracy test is the linear advection-reaction system from \cite{HundsdorferRuuth2007}. The equations are
\begin{eqnarray}
\partial_t u + \alpha_1\,\partial_x u &=& -k_1 u + k_2 v + s_1\,,\\
\partial_t v + \alpha_2\,\partial_x v &=& k_1 u - k_2 v + s_2
\end{eqnarray}
for $0<x<1$ and $0<t\le 1$, with parameters
\[
\alpha_1=1,\;\alpha_2=0,\;k_1=10^6,\;k_2=2k_1,\;s_1=0,\;s_2=1,
\]
and with the following initial and boundary conditions:
\[
u(x,0)=1+s_2x,\;v(x,0)=\frac{k_1}{k_2}u(x,0)+\frac{1}{k_2}s_2,\;
u(0,t)=1-\sin(12t)^4\,.
\]
Note that there are no boundary conditions for $v$ since $\alpha_2$
is set to be zero.

Fourth-order finite differences on a uniform mesh consisting of $m=400$
nodes are applied in the interior of the domain. At the boundary we can
take third-order upwind biased finite differences, which here does not
affect an overall accuracy of four \cite{HundsdorferRuuth2007}
and gives rise for a spatial error of $1.5\cdot 10^{-5}$.

In the IMEX setting, the reaction is treated implicitly and all other
terms explicitly. In order to guarantee that errors of the initial
values do not affect the computations, accurate initial values are computed
by {\sc ode15s} from {\sc Matlab} with high tolerances. We have used step sizes
$\st=(4,2,1,0.5,0.25,0.1,0.05,0.025)\cdot 10^{-3}$, and compared the
numerical values at the final time $T=1$ with an accurate reference
solution in the $l_2$-vector norm ($\|v\|^2=\sum_iv_i^2$) as used in
\cite{LangHundsdorfer2017}.
The results are plotted in Figure~\ref{fig:AdvReac}.

\begin{figure}[ht]
\centering
\includegraphics[width=0.8\textwidth]{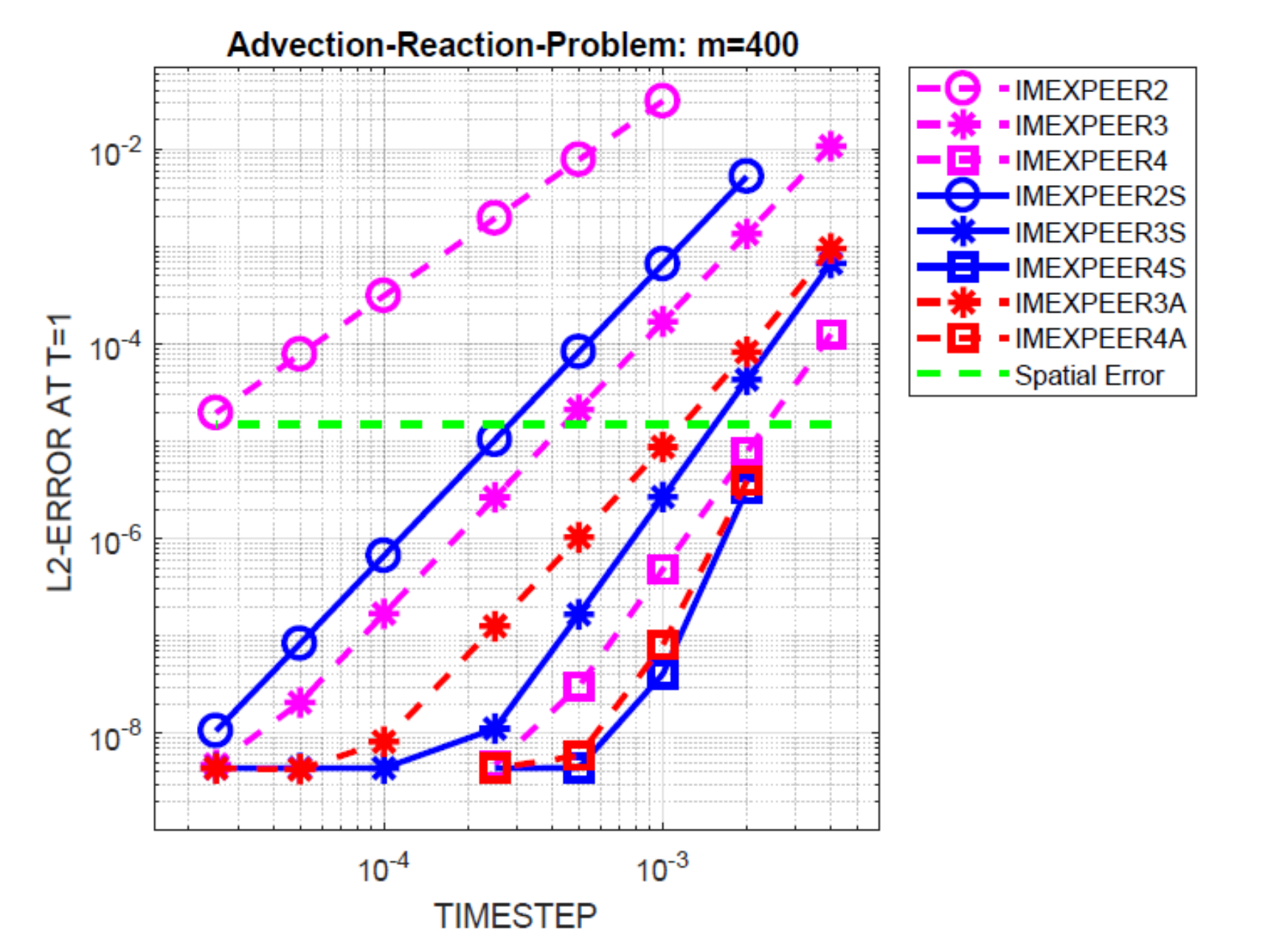}
\parbox{13cm}{
\caption{
Advection-Reaction-Problem: $l_2$-errors at $T=1$
of the total concentration vs. time step sizes, $m=400$.
Comparison of IMEX-Peer methods. To
keep clarity for the other methods, results for the $4$-stage
methods applied with $\st\le 10^{-4}$ are not shown. The corresponding errors
stagnate at $4\cdot 10^{-9}$. Further missing results correspond to
failures or unstable numerical solutions.}
\label{fig:AdvReac}
}
\end{figure}

All methods show their theoretical orders, but the larger error
constants of the super-convergent IMEX-Peer schemes compared
to the IMEX-Peer(s) family are again apparent. The $3$-stage
methods and IMEX-Peer4 still deliver satisfactory results for the
largest time step, $\st=4\cdot 10^{-3}$, whereas the others fail.
The similar asymptotic behaviour of the higher-order methods
shows order reduction for smaller time steps,
which was also observed in \cite{HundsdorferRuuth2007} as an inherent
issue for very high-accuracy computations. However, this effect
appears on a level far below the spatial discretization error.

\subsection{Nonlinear Two-Dimensional Gravity Waves}
This problem is taken from \cite{SoleimaniKnothWeiner2017}.
Let $u(x,z)$ be the velocity of the gravity waves with $x$ and $z$
denoting the horizontal and vertical coordinates in $\si{\metre}$, respectively.
The gravity waves are generated by a localized region of a non-divergent
forcing in a stratified shear flow. The horizontal background wind
(in \si{\metre\per\second}) is given by
\[ u_0(z ) = 5 + \frac{z}{1000} + 0.4\left( 5-\frac{z}{1000}\right)\left( 5+\frac{z}{1000}\right)\]
and the waves are forced by the curl of a non-divergent stream function
(in \si{\metre\squared\per\second})
\[ \psi(x,z,t) = \psi_0\left(\frac{\pi x}{L_x}\right)\sin(\omega t)
\exp\left[ -\left(\frac{\pi x}{L_x}\right)^2-
\left(\frac{\pi z}{L_z}\right)^2\right].\]
The parameters used here are $\psi_0=80\;\si{\metre\squared\per\second}$,
$L_x=10\cdot10^3\,\si{\metre}$ and $L_z=2.5\cdot10^3\,\si{\metre}$.
The governing system of equations reads
\begin{align}
\frac{Du}{Dt} + \frac{\partial P}{\partial x} & =
  - \frac{\partial \psi}{\partial z} + \frac{u_0(z)- \bar{u}(z,t)}{\tau}\,, \\[1mm]
\frac{Dw}{Dt} + \frac{\partial P}{\partial z} & =
  b+ \frac{\partial \psi}{\partial x}\,,\\[1mm]
\frac{Db}{Dt} + N^2w & = 0,\\[1mm]
\frac{DP}{Dt} + c_s^2\left( \frac{\partial u}{\partial x} +
\frac{\partial w}{\partial z}\right) & = 0\,,
\end{align}
with horizontally averaged mean flow $\bar{u}$,
$x\in [-15\,L_x,15\,L_x]$, $z\in [-2\,L_z,2\,L_z]$
and $t\in [0\,\si{\second},3000\,\si{\second}]$. The constants are
\[\omega = 0.005\,\si{\per\second}, \quad N = 0.02
\quad \textrm{and} \quad c_s = 350\,\si{\metre\per\second}\quad \textrm{(speed of sound)}.\]
\begin{figure}[t]
\centering
\includegraphics[width=0.8\textwidth]{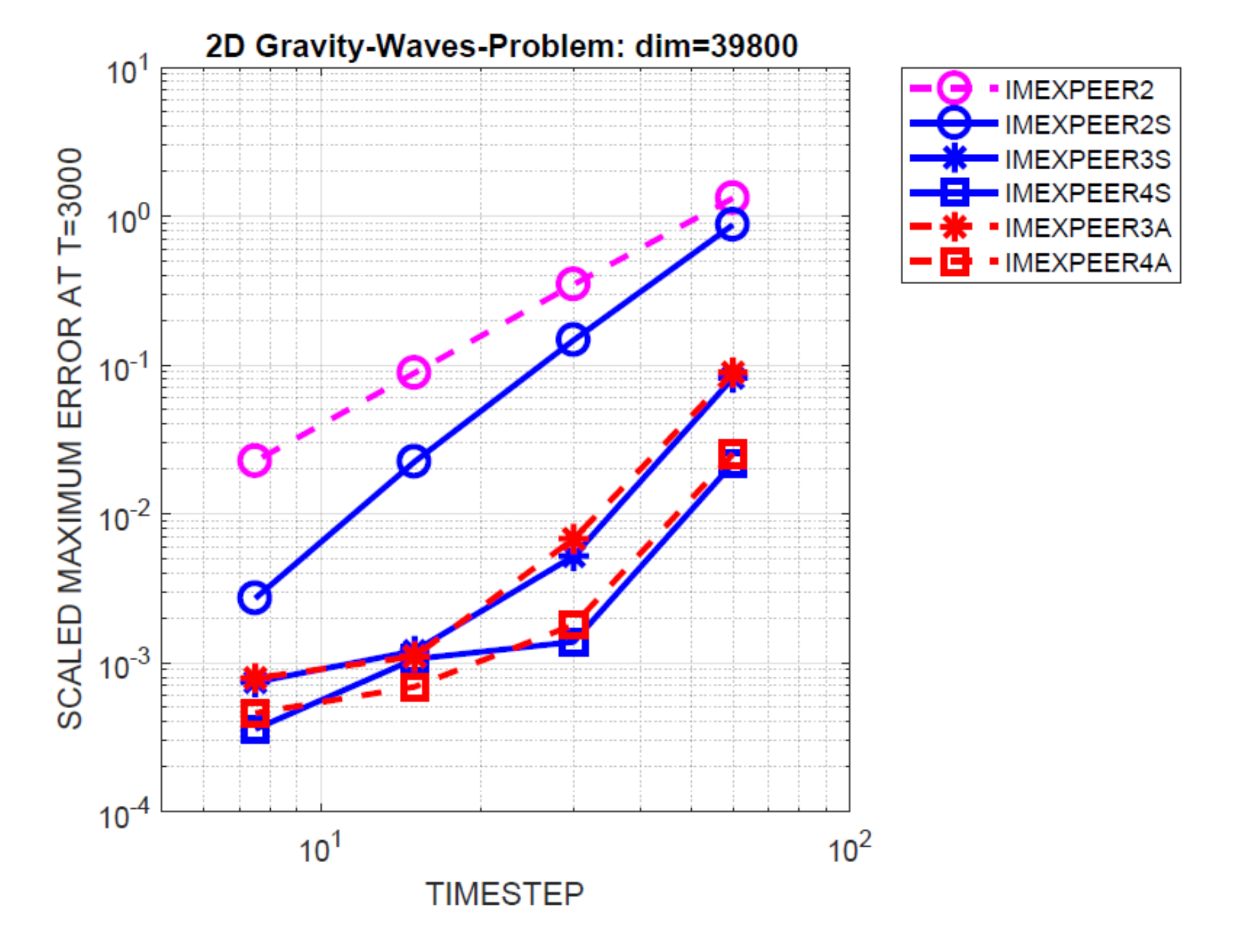}
\parbox{13cm}{
\caption{\small
Nonlinear Two-Dimensional Gravity Waves: Scaled maximum errors at $T=3000\,\si{\second}$
taken over all system components vs. time step sizes. Comparison of IMEX-Peer methods.}
\label{fig:Gravity}
}
\end{figure}
For the spatial discretization, we set
\begin{align*}
\delta_{nx}f(x) & = \frac{f(x + n\Delta x/2) - f(x - n\Delta x/2)}{n\Delta x}\,,\\[2mm]
\langle f(x)\rangle^{nx} & = \frac{f(x + n\Delta x/2) + f(x - n\Delta x/2)}{2}\,,\\[2mm]
\psi_x & = \frac{\psi(x+\Delta x,z,t)-\psi(x,z,t)}{\Delta x}\,,\\[2mm]
\psi_z & = \frac{\psi(x,z+\Delta z,t)-\psi(x,z,t)}{\Delta z}
\end{align*}
and obtain the semi-discretized system of ODEs
\begin{align}
\label{SpatGravityFirst}
\frac{\partial u}{\partial t}  + \frac1 2 \delta_{2x}(u^2) + \langle\langle w\rangle^x\delta_z u\rangle^z - \frac{u_0(z)-\bar{u}(z,t)}{\tau} + \psi_z + K \Delta_h u & = -\delta_x P \\[1mm]
\frac{\partial w}{\partial t} + \langle\langle u\rangle^z \delta_x w\rangle^x + \frac{1}{2} \delta_{2z}(w^2) - \psi_x + K\Delta_h w & = - \delta_z P + b \\[1mm]
\frac{\partial b}{\partial t} + \langle\langle u\rangle^z \delta_x b\rangle^x + \langle\langle w\rangle^z \delta_z b\rangle^z  + K \Delta_h b & = - N^2w \\[1mm]
\label{SpatGravityLast}
\frac{\partial P}{\partial t} + \langle u \delta_x P\rangle^x + \langle w \delta_z P\rangle^z & = -c_s^2(\delta_x u + \delta_z w),
\end{align}
where $K = 4.69 \cdot 10^{-4}\,\si{\per\second}$. We add a fourth-derivative hyper-diffusion term
\[ \Delta_h : = \left( (\Delta x\delta_x)^2 +(\Delta z \delta_z)^2\right)^2\]
to prevent nonlinear instability and to impose a simple parametrization of turbulent mixing.
In our calculation, we have used $m_x=100$ and $m_z=100 $ nodes for the spatial grid. This yields a system of dimension $39800$. A reference solution is computed by Matlab's {\sc ode15s} with high precision.
For more details, we refer to \cite{SoleimaniKnothWeiner2017}, see also \cite{DuranBlossey2012}.

We identify the right hand side of \eqref{SpatGravityFirst}--\eqref{SpatGravityLast} as the stiff part of the ODE system, which is treated implicitly.  The eigenvalues of its Jacobian
are all complex and lie approximately in $[-7\i,7\i]$. This especially shows that $A$-stability is necessary for a stable computation of the gravity waves. It also explains why the methods
IMEX-Peer3 and IMEX-Peer4 \cite{LangHundsdorfer2017} with $\alpha<90^\circ$, i.e., the implicit
method is not A-stable, fail for this problem and are not shown in Figure~\ref{fig:Gravity}.
We have used step sizes $\st=60,30,15,7.5\,\si{\second}$ and compared the numerical values at the final time $T=3000\,\si{\second}$ with an accurate reference solution $y_{ref}$ in the scaled maximum norm,
i.e., we set $err=\max_{i=1,\ldots,n}|y_i-y_{ref,i}|/(1+|y_{ref,i}|)$, see Figure~\ref{fig:Gravity}.

The super-convergent IMEX-Peer2s clearly shows its additional order.
IMEX-Peer3s and (the nearly super-convergent) IMEX-Peer3a deliver almost identical results due to their similar
error constants, whereas super-convergence can only be observed for larger time steps. Both
4-stage IMEX-Peer methods show an unpredictable behaviour, but give the smallest errors.
Similar results are presented in \cite{SoleimaniKnothWeiner2017}.

\section{Conclusion}
We have developed a new family of $s$-stage IMEX-Peer methods,
applying the idea of extrapolation from \cite{CardoneJackiewiczSanduZhang2014,LangHundsdorfer2017} to a broader
class of implicit Peer methods which also include function values from
the previous step. These additional degrees of freedom allow the
construction of super-convergent IMEX-Peer methods of order $s\!+\!1$
with A-stable implicit part. A-stability is important to solve problems
with function contributions that have large imaginary eigenvalues in
the spectrum of their Jacobians. We analysed the property of
super-convergence and gave sufficient conditions for it. These
conditions are motivated by the definition of the order of consistency
for general linear methods and are more general as those commonly
used in Peer literature. Linear stability properties were carefully
examined to derive new super-convergent $s$-stage IMEX-Peer methods
with order $p\!=\!s\!+\!1$ for $s\!=\!2,3,4$. We employed
the {\sc Matlab}-routine {\it fminsearch} with varying objective
functions and starting values to find suitable methods with sufficiently
large stability regions and small error constants. However,
the new properties, super-convergence and A-stability, go along with
significantly smaller stability regions and larger error constants,
compared to the IMEX-Peer methods from \cite{LangHundsdorfer2017}.

In a detailed comparison with recently proposed IMEX-Peer methods
in \cite{LangHundsdorfer2017,SoleimaniKnothWeiner2017}, the new
methods showed their super-convergence property and
performed better in many cases. They gave equally well results for
a two-dimensional gravity wave problem used in
\cite{SoleimaniKnothWeiner2017} to demonstrate the importance
of A-stability of the implicit part, whereas the higher-order methods from
\cite{LangHundsdorfer2017} failed due to their lack of this property.
In future work, we will extend our construction principles to
variable step size sequences with local error control. This is
especially important when it comes to solve real-life applications
with often largely varying time scales, see e.g.
\cite{GerischLangPodhaiskyWeiner2009,GottermeierLang2009} for
linearly implicit Peer methods with variable step sizes.

\section{Acknowledgement}
The authors would like to thank B.~Soleimani, O.~Knoth and R.~Weiner
for making the nonlinear two-dimensional gravity wave problem used in
\cite{SoleimaniKnothWeiner2017}
available for our comparison. J.~Lang was supported by the
German Research Foundation within the collaborative research center
TRR154 ``Mathematical Modeling, Simulation and Optimisation Using
the Example of Gas Networks'' (DFG-SFB TRR154/1-2014, TP B01),
the Graduate School of Excellence Computational Engineering (DFG GSC233),
and the Graduate School of Excellence Energy Science and
Engineering (DFG GSC1070). We thank Willem for his long-time warm friendship
and inspiring scientific collaboration. He will be greatly missed.

\bibliographystyle{plain}
\bibliography{bibimexpeersuper}

\end{document}